\documentclass[11pt]{article}
\usepackage{latexsym}
\newcounter{dim}
\title{}
\author{}
\date{}
\begin{document}
\voffset=-1.3truein
\hoffset=-.5truein
\def\skp#1{\vskip#1cm\relax}
\def\sumk#1#2{\sum_{k = #1}^{\infty} #2}
\def\sumn#1#2{\sum_{n = #1}^{\infty} #2}
\def\tn{{\bf T$^{n}$}}
\def\pn{$P^{n}\/$}
\def\ml{M^{2n}(\lambda)}
\def\cn{{\bf C$^{n}$}}
\def\s1{{\bf S$^{1}$}}
\def\a1{${\cal{A}}(1)$}
\def\ld{\lambda}
\def\z2{{\bf Z}$_{2}$}
\def\zl2{{\bf Z}$_{(2)}$}
\def\block{\rule{2.4mm}{2.4mm}}
\def\rplus{{\bf R$_{+}}}
\def\nd{\noindent}
\newtheorem{theorem}{Theorem}
\maketitle
\skp{.1}
\centerline{{{\large {\sc THE $KO$-THEORY OF TORIC MANIFOLDS}}}}
\skp{.65}
\centerline{{{\small {\sc ANTHONY BAHRI AND MARTIN BENDERSKY}}}}
\skp{.45}
\centerline{\today}
\skp{0.9}
\begin{abstract}
\nd Toric manifolds, a topological generalization of smooth projective
toric varieties, are
determined by an $n$-dimensional simple convex polytope and a function
from the set of codimension-one faces into the primitive vectors of
an integer lattice. Their cohomology was determined by Davis and
Januszkiewicz in 1991 and corresponds with the theorem of Danilov-
Jurkiewicz in the toric variety case. Recently it has been shown by
Buchstaber and Ray that they generate the complex cobordism ring.
We use the Adams spectral sequence to compute the $KO$-theory of 
all toric manifolds and certain singular toric varieties. 
\end{abstract}
\skp{.3}
\centerline{{\bf 1. Introduction}}
\skp{0.1}
%
%
We take as our definition of toric manifold the construction of Davis
and Januszkiewicz (\cite{dj}, section $1.5$). 
%
%
Let \pn $\;$ be an $n$-dimensional, simple (at each vertex, $n$ 
codimension-one faces meet), convex polytope. 
Set
$${\cal F} = \{F_1,F_2,\ldots,F_m\}$$
the set of codimension-one faces of \pn. The fact that \pn $\;$ is 
simple implies that every codimension-$l$ face $F$ can be 
written uniquely as
$$F = F_{i_{1}} \cap F_{i_{2}} \cap \cdots \cap F_{i_{l}}$$
where the $F_{i_{j}}$ are codimension-one faces containing $F$. Let
$$\lambda : {\cal F} \rightarrow {\mbox{{\bf Z}}}^n$$
be a function into an $n$-dimensional integer lattice satisfying the 
condition that whenever 
$F = F_{i_{1}} \cap F_{i_{2}} \cap \cdots \cap F_{i_{l}}$ then 
$\lambda(F_{i_{1}}),\lambda(F_{i_{2}}), \ldots ,\lambda(F_{i_{l}})$ span 
an $l$-dimensional submodule of {\bf Z}$^n$ which is a direct summand.
Next, regarding {\bf R}$^n$ as the Lie algebra of \tn, we see that
$\lambda$ associates to each codimension-$l$ face $F$ of \pn $\;$ a rank-$l$
subgroup $G_F \subset$ \tn. Finally, let $p \in$ \pn $\;$ and $F(p)$ be 
the unique face with $p$ in its relative interior. Define an equivalence
relation $\sim$ on \tn $\times$ \pn $\;$ by $(g,p) \sim (h,q)$ if and only
if $p = q$ and $g^{-1}h \in G_{F(p)} \cong$ {\bf T}$^l$. Set
$$M^{2n}(\lambda) = {\mbox{{\bf T}}^n \times P^n}\big/\!\sim$$ 
$\ml$ is a smooth, closed, connected, $2n$-dimensional manifold with a 
\tn action induced by left translation (\cite{dj}, page 423). There is 
a projection
$$\pi : \ml \rightarrow P^n$$
induced from the projection \tn $\times$ \pn $\rightarrow$ \pn.

Following \cite{dj}, we note that every toric manifold has this 
description, in particular, every smooth projective toric variety 
does too. Recently, Buchstaber and Ray \cite{br} have shown that toric 
manifolds generate the complex cobordism ring.

Here is a simple example selected from the list in \cite{dj}.
Let $n = 2$ and $P^2$ be a square. Here
${\cal F} = \{F_1,F_2,F_3,F_4\}$ consists of four codimension-one 
faces. Define $\lambda : {\cal F} \rightarrow {\mbox{{\bf Z}}}^2$ 
as in the diagram below.
\skp{0.1}
\begin{picture}(14,4)(0,0)
\put(3,1){\line(1,0){1.5}}
\put(3,2.5){\line(1,0){1.5}}
\put(3,1){\line(0,1){1.5}}
\put(4.5,1){\line(0,1){1.5}}
\put(3,1){\circle*{0.09}}
\put(3,2.5){\circle*{0.09}}
\put(4.5,2.5){\circle*{0.09}}
\put(4.5,1){\circle*{0.09}}
\put(1,1.7){\makebox(0,0)[bl]{{\mbox{\scriptsize $\lambda(F_1) = (0,1)$}}}}
\put(4.7,1.7){\makebox(0,0)[bl]{{\mbox{\scriptsize $\lambda(F_3) = (-1,1)$}}}}
\put(3.75,.6){\makebox(0,0){{\mbox{\scriptsize $\lambda(F_2) = (1,0)$}}}}
\put(3.75,2.8){\makebox(0,0){{\mbox{\scriptsize $\lambda(F_4) = (1,-2)$}}}}
\put(7.5,1.7){\makebox(0,0)[bl]{yields}}
\put(9.2,1.7){\makebox(0,0)[bl]{$M^{4}(\lambda) \cong CP^2 \# CP^2$}}
\end{picture}
\skp{0.1}
\nd Davis and Januszkiewicz point out that $CP^2 \# CP^2$ 
is a toric manifold but does not have an almost complex 
structure and so cannot be a toric variety. Our main results are:
\skp{.4}
\nd {\bf Theorem 1.} {\em The Adams spectral sequence for the
real connective\/} $KO$-{\em theory of the toric manifold\/}
$M^{2n}(\lambda)$ {\em collapses\/}.
\skp{0.2}
\nd {\bf Corollary 2.} $KO^{*}{M^{2n}(\lambda)}$ {\em is determined
by the\/} $\bmod{\,2}$ {\em cohomology ring of\/} $M^{2n}(\lambda)$.
{\em In particular, the\/} $KO$-{\em theory depends only the values
of\/} $\lambda \bmod{\,2}.$
\skp{0.2}
\nd Our methods yield the additional result that the theorem remains
true for certain singular toric varieties, of real dimension less 
than 12.
\skp{0.2}
\nd We note that the $K$-theory of toric varieties has been computed
by Robert Morelli in \cite{rm} 
\skp{0.4}
\nd {\bf Acknowledgement.\/} We are grateful to Ciprian Borcea for
his encouragement and helpful comments and for introducing us to this 
subject through a series of fine seminars he gave on toric varieties at 
Rider University. We would like also to thank Bob Bruner for several 
helpful conversations.
\skp{.4}

%

\skp{.4}
\centerline{{\bf 2. Homology and Cohomology of}
$\displaystyle M^{2n}(\lambda)$}
\skp{0.1}
\nd In order to compute the $KO$-theory of $M^{2n}(\lambda)$ we shall 
need the computation of its homology from \cite{dj}. To state their result
we recall certain numbers defined in terms of the combinatorics of \pn.
Let $f_{i}$ be the number of faces of \pn $\;$ of codimension $(i + 1)$. Define
numbers $h_{i}$ by the equality of polynomials in $t$
$$(t - 1)^{n} + \sum_{i = 0}^{n - 1}{f_{i}(t - 1)^{n - 1 - i}}
= \sum_{i = 0}^{n}{h_{i}t^{n - i}}$$
$(h_{0},\ldots,h_{n})$ is called the {\it h-vector\/} of \pn. Notice
$h_0 = h_n = 1$ and 
$$\sum_{i = 0}^{n}{h_{i}} = f_{n - 1} = \mbox{the number of vertices of}
\; P^{n}$$
For each $k$-face $F$ of \pn $\;$ we have a connected $2k$-dimensional 
submanifold $M_F$ of $\ml$ defined by $M_F = \pi^{-1}(F)$.
\skp{0.4}
\nd {\bf Theorem 3. [M. Davis and T. Januszkiewicz \cite{dj}]}
{\em The group\/} $H_{*}(\ml$;{\bf Z}) {\em is independent of the 
function\/} $\lambda$. {\em Specifically\/},
$$\begin{array}{ccl}
H_{2i + 1}(\ml; \mbox{\bf Z}) & = & 0\\
\\
H_{2i}(\ml; \mbox{\bf Z})     & = & \mbox{free of rank} \; h_{i}
\end{array}$$
The group $H^{2l}(\ml; \mbox{\bf Z})$ {\em is generated by the 
Poincar\'{e} duals
of classes of the form\/} $[M_F]$ {\em with \/} $F$ {\em a face of 
codimension} $l$.
{\em As a ring\/}, $H^{*}(\ml; \mbox{\bf Z})$ {\em is generated by the 
degree-two classes dual to\/} $[M_F]$ {\em with\/} $F$ {\em a 
face of codimension one\/}. \block 
\skp{0.3}
\nd The ring structure of $H^{*}(\ml; \mbox{\bf Z})$ is determined from the
Serre spectral sequence of the fibration
$$\ml \rightarrow BP^{n} \rightarrow B{\mbox{\bf T}}^{n}$$ where $BP^{n}$
denotes the Borel construction
$$BP^{n} = E{\mbox{\bf T}}^{n} \times_{{\mbox{\bf T}}^{n}} \ml$$
Let $v_1,v_2,\ldots,v_m$ denote the degree-two generators of 
$H^{*}(\ml; \mbox{\bf Z})$, one for each codimension-one face of \pn. We
need to define two ideals of relations in {\bf Z}$[v_1,v_2,\ldots,v_m]$,
$I$ and $J$. 

Let $K$ be the simplicial complex {\em dual\/}
to \pn. That is, an $(n - 1)$-dimensional simplicial complex with vertex 
set $\cal{F}$, the set of codimension-one faces of \pn. A set of $(k + 1)$
elements in $\cal{F}$, $\{F_{i_{0}},\ldots,F_{i_{k}}\}$ span a $k$-simplex
in $K$ if and only if $F_{i_{0}}\cap \cdots \cap F_{i_{k}} \neq \phi$. 
The ideal $I$ is the homogenous ideal of relations generated 
by all square free 
monomials of the form $v_{i_{1}}\cdots v_{i_{s}}$, where 
$\{v_{i_{1}}, \ldots, v_{i_{s}}\}$ does not span a simplex in $K$.

The ideal $J$ is defined in terms of the function $\lambda$. Let 
$\{e_1,\ldots, e_m\}$ be the standard basis of {\bf Z}$^m$. Then, 
identifying the codimension-one face $F_{i}$ with $e_i$, we can regard
$$\lambda : {\cal F} \rightarrow {\mbox{{\bf Z}}}^n$$
as a linear map {\bf Z}$^{m}$ $\rightarrow$ {\bf Z}$^{n}$ given by an
$m \times n$ matrix $(\lambda_{ij})$.
In example 3 above, the linear map 
$\lambda :$ {\bf Z}$^{4}$ $\rightarrow$ {\bf Z}$^{2}$
is the matrix
$$\lambda = \pmatrix{0&1&-1&1\cr
                     1&0&1&-2\cr}$$
The ideal of relations $J$ is determined by the system of equations
$$\begin{array}{ccc}
\ld_{11}v_{1} + \ld_{12}v_{2} + \ldots + \ld_{1m}v_{m} & = & 0\\
\ld_{21}v_{1} + \ld_{22}v_{2} + \ldots + \ld_{2m}v_{m} & = & 0\\
\vdots                                       & \vdots & \vdots\\
\ld_{ni}v_{1} + \ld_{n2}v_{2} + \ldots + \ld_{nm}v_{m} & = & 0
\end{array}$$                
\skp{0.4}
\nd {\bf Theorem 4. [M. Davis and T. Januszkiewicz \cite{dj}]}
{\em As rings\/}
$$H^{*}(\ml;\mbox{{\bf Z}}) = 
{\mbox{{\bf Z}}[v_1,v_2,\ldots,v_m]} \bigg/(I + J) \quad \mbox{\block}$$
\skp{0.3}
\newpage
\nd As an illustration, we compute $H^{*}(M^{4}(\ld);\mbox{\bf Z})$ with 
$M^{4}(\lambda) \cong CP^2 \# CP^2$, the example from the introduction.
The dual of $P^2$ is a one-dimensional simplicial complex $K$ with
vertices $\{v_1,v_2,v_3,v_4\}$.
\skp{0.3}
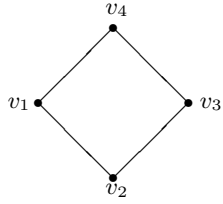
\begin{picture}(14,3)(0,0)
\put(1.5,1.5){\line(1,1){1}}
\put(1.5,1.5){\line(1,-1){1}}
\put(2.5,2.5){\line(1,-1){1}}
\put(2.5,.5){\line(1,1){1}}
\put(1.5,1.5){\circle*{0.09}}
\put(2.5,2.5){\circle*{0.09}}
\put(3.5,1.5){\circle*{0.09}}
\put(2.5,.5){\circle*{0.09}}
\put(1.1,1.4){\makebox(0,0)[bl]{{\mbox{\footnotesize $v_1$}}}}
\put(3.65,1.4){\makebox(0,0)[bl]{{\mbox{\footnotesize $v_3$}}}}
\put(2.4,.25){\makebox(0,0)[bl]{{\mbox{\footnotesize $v_2$}}}}
\put(2.4,2.65){\makebox(0,0)[bl]{{\mbox{\footnotesize $v_4$}}}}
\put(5,2){\makebox(0,0)[bl]{{\mbox{\small $\{v_1,v_3\}$
                                   does not span a simplex}}}}
\put(5,1.5){\makebox(0,0)[bl]{{\mbox{\small $\{v_2,v_4\}$
                                   does not span a simplex}}}}
\put(5,.9){\makebox(0,0)[bl]{{\mbox{\small so 
  $I = \;<v_{1}v_{3}, v_{2}v_{4}>\; 
       \subset$ {\bf Z}$[v_1,v_2,v_3,v_4]$}}}}
\end{picture}
\skp{0.2}
\nd The relations $J$ are read off from the matrix $\ld$ above
\skp{0.1}
$$\left. \begin{array}{ccc}
v_2 - v_3 + v_4  & = & 0\\
v_1 + v_3 - 2v_4 & = & 0
\end{array}\right\} \Rightarrow 
\begin{array}{ccc}
v_3 & = & v_2 + v_4\\
v_1 & = & v_4 - v_2
\end{array}$$
\skp{0.2}
\nd Choosing generators $v_2, v_4 \in H^{2}(M^{4}(\ld);\mbox{\bf Z})$ 
we get
\skp{0.1}
$$\begin{array}{cccl}
H^{0}(M^{4}(\ld);\mbox{\bf Z})  & = & \mbox{\bf Z} &\\
H^{2}(M^{4}(\ld);\mbox{\bf Z})  & = & \mbox{\bf Z} \oplus \mbox{\bf Z} &
                       <v_2, v_4>\\
H^{4}(M^{4}(\ld);\mbox{\bf Z})  & = & \mbox{\bf Z}                     &
                       <v_{2}^{2} = v_{4}^{2}>\\
H^{i}(M^{4}(\ld);\mbox{\bf Z})  & = & 0                                &
                       i > 4, \quad v_{{i}_{1}}v_{{i}_{2}}v_{{i}_{3}} = 0,
                       \quad i_j \in \{2,4\}
\end{array}$$                       
\skp{.4}
\centerline{{\bf 3. The Action of the Steenrod Algebra}}
\skp{0.1}
\nd For our calculation, we require the structure of 
$H^{*}(M^{2n}(\ld); \mbox{\bf Z$_{2}$})$ as a module over the subalgebra
${\cal{A}}(1)$ , generated by $Sq^1$ and $Sq^2$, of the $\bmod{\,2}$ Steenrod 
algebra $\cal{A}$. Let $S^{0}$ denote the \a1 module consisting of a 
single class in dimension 0 and the trivial action of $Sq^1$ and $Sq^2$.
Denote by ${\cal M}$ the \a1 module with a class $x$ in dimension 0, a class
$y$ in dimension 2 and the action given by $Sq^2(x) = y$.
\skp{0.4}
\nd {\bf Lemma 5.} {\em Let\/} $X$ {\em be a space with\/} 
$H^*(X;\mbox{\bf Z$_2$})$
{\em concentrated in even degrees. Then, as an\/} \a1
{\em module, $H^*(X;\mbox{\bf Z$_2$})$ is isomorphic to a direct sum of 
suspended copies of\/}
$S^0$ {\em and\/} ${\cal M}$. {\em Furthermore, the splitting is natural with 
respect to maps of spaces\/}.
\skp{0.25}
\nd {\bf Proof:\/} The sequence
$$\to H^{2n - 2}(X;\mbox{\bf Z$_2$}) \quad 
{\mathop {\longrightarrow} \limits^{Sq^2}} 
\quad H^{2n}(X;\mbox{\bf Z$_2$}) \to $$
is a chain complex since $Sq^{2}Sq^{2} = Sq^{3}Sq^{1} = 0$ because
$H^*(X;\mbox{\bf Z$_2$})$ is concentrated in even degrees.
Its homology is defined to be the ``$Sq^2$ homology of $X$" and 
is denoted $$H_*(X ; Sq^2).$$
Let $A_{2n}=$ Ker$\{ Sq^2 : H^{2n}(X;\mbox{\bf Z$_2$}) \to 
H^{2n + 2}(X;\mbox{\bf Z$_2$}) \}$.  Then 
$H^{2n}(X;\mbox{\bf Z$_2$}) \approx A_{2n} \oplus B_{2n}$ for some 
vector subspace $B_{2n}$.
Define $C_{2n} \subseteq A_{2n}$ to be 
Im$\{ Sq^2: H^{2n - 2}(X;\mbox{\bf Z$_2$}) \to 
H^{2n}(X;\mbox{\bf Z$_2$}) \}$. 
Then $A_{2n} \approx C_{2n} \oplus D_{2n}$ for some vector subspace 
$D_{2n}$.  Hence we have $H^{2n}(X;\mbox{\bf Z$_2$}) \approx 
C_{2n} \oplus D_{2n} \oplus 
B_{2n}$ with  $H_{2n}(X ;Sq^2) \approx D_{2n}$ and 
$Sq^2 : B_{2n-2} \to C_{2n}$ an isomorphism.  
The lemma now follows since $D_{2n}$ generates copies of suspensions 
of $S^0$ 
and $B_{2n} (\approx C_{2n + 2})$ generates suspensions of ${\cal M}$. 
The naturality follows since 
$H_*(X; Sq^2)$ and $C_*$ are natural. \block
\skp{0.27}
\newpage
An algorithm allows us to determine the \a1 module structure of
$H^*(X;\mbox{\bf Z$_2$})$ explicitly.
Let $\{u_{(2,1)},u_{(2,2)}, \ldots, u_{(2,s_{2})}\}$ be a \z2
basis for $H^{2}(X;\mbox{\bf Z$_2$})$. We construct a new basis 
$\{w_{(2,1)},w_{(2,2)}, \ldots, w_{(2,s_{2})}\}$ which will yield the
decomposition above. Set $w_{(2,1)} = u_{(2,1)}$.
If $Sq^{2}u_{(2,2)} = Sq^{2}w_{(2,1)}$ set 
$w_{(2,2)} = w_{(2,1)} + u_{(2,2)}$, else $w_{(2,2)} = u_{(2,2)}$. 
Suppose now that $w_{(2,t - 1)}$
has been defined. If $Sq^{2}u_{(2,t)}$ is linearly independent of
$\{Sq^{2}w_{(2,1)},Sq^{2}w_{(2,2)},\ldots ,Sq^{2}w_{(2,t - 1)}\}$ set
$w_{(2,t)} = u_{(2,t)}$. Otherwise, if
$$Sq^{2}u_{(2,t)} = Sq^{2}w_{(2,i_{1})} + Sq^{2}w_{(2,i_{2})} + \ldots +
Sq^{2}w_{(2,i_{t})}$$
set $w_{(2,t)} = u_{(2,t)} + w_{(2,i_{1})} + \ldots + w_{(2,i_{t})}$.
Next, reorder the set $\{w_{(2,1)},w_{(2,2)}, \ldots, w_{(2,s_{2})}\}$
so that $Sq^{2}w_{(2,j)} = 0$ for $j = 1, \ldots, t_2$ and 
$Sq^{2}w_{(2,j)} \neq 0$ for $j = t_{2} + 1, \ldots, s_2$. Set
$d_{(2,j)} = w_{(2,j)}$ for $j = 1, \ldots, t_2$ and 
$b_{(2,j)} = w_{(2,t_{2} + j)}$ for $j = 1, \ldots, s_2 - t_2$. 
So, in the notation above, 
$$D_2 = \{d_{(2,1)},d_{(2,2)}, \ldots, d_{(2,t_2)}\}$$
and
$$B_2 = \{b_{(2,1)},b_{(2,2)}, \ldots, b_{(2,s_2 - t_2)}\}$$
Of course, $C_2 = \phi$ and $C_4 \approx B_2$.
Now suppose that $A_{2n - 2}, B_{2n - 2}$ and $C_{2n - 2}$ have
been constructed. Set 
$$C_{2n} = \{Sq^{2}b_{(2n - 2,1)},Sq^{2}b_{(2n - 2,2)}, \ldots, 
Sq^{2}b_{(2n - 2,s_{2n - 2} - t_{2n - 2})}\} \approx B_{2n -2}.$$ 
The elements of
$C_{2n}$ are linearly independent by the construction of $B_{2n - 2}$.
Choose {\em any\/} extension of $C_{2n}$ to a basis of 
$N^{2n} = H^{2n}(X;\mbox{\bf Z$_2$})$. Denote the basis by
$$C_{2n} \cup \{u_{(2n,1)},u_{(2n,2)}, \ldots, u_{(2n,s_{2n})}\}$$
Finally, repeat the process above on the set
$$\{u_{(2n,1)},u_{(2n,2)}, \ldots, u_{(2n,s_{2n})}\}$$
to produce $B_{2n}$ and $D_{2n}$.
\skp{0.3}
\nd Diagrammatically, the \a1 module structure looks like
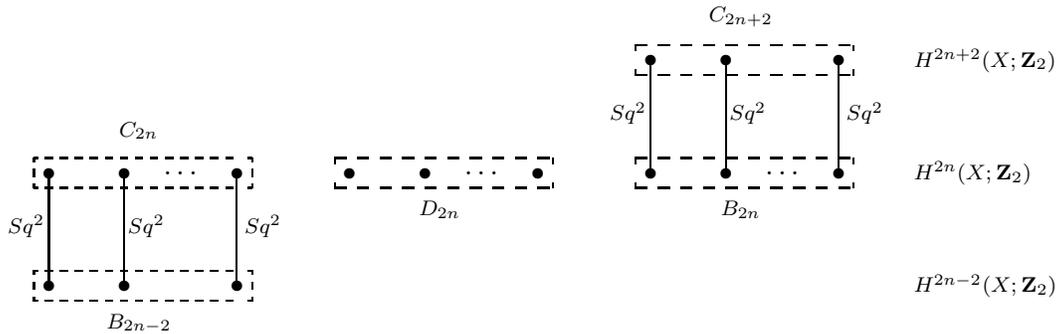
\skp{0.3}
\begin{picture}(14,6)(0,0)
\put(.5,3){\circle*{.15}}
\put(1.5,3){\circle*{.15}}
\put(3,3){\circle*{.15}}
\put(.5,1.5){\circle*{.15}}
\put(1.5,1.5){\circle*{.15}}
\put(3,1.5){\circle*{.15}}
\put(4.5,3){\circle*{.15}}
\put(5.5,3){\circle*{.15}}
\put(7,3){\circle*{.15}}
\put(8.5,3){\circle*{.15}}
\put(9.5,3){\circle*{.15}}
\put(11,3){\circle*{.15}}
\put(8.5,4.5){\circle*{.15}}
\put(9.5,4.5){\circle*{.15}}
\put(11,4.5){\circle*{.15}}
\put(.5,1.5){\line(0,1){1.5}}
\put(1.5,1.5){\line(0,1){1.5}}
\put(3,1.5){\line(0,1){1.5}}
\put(8.5,3){\line(0,1){1.5}}
\put(9.5,3){\line(0,1){1.5}}
\put(11,3){\line(0,1){1.5}}
\put(2.25,3){\makebox(0,0){\mbox{$\ldots$}}}
\put(6.25,3){\makebox(0,0){\mbox{$\ldots$}}}
\put(10.25,3){\makebox(0,0){\mbox{$\ldots$}}}
\put(.2,2.15){\makebox(0,0)[b]{{\mbox{\scriptsize $Sq^2$}}}}
\put(1.55,2.15){\makebox(0,0)[bl]{{\mbox{\scriptsize $Sq^2$}}}}
\put(3.1,2.15){\makebox(0,0)[bl]{{\mbox{\scriptsize $Sq^2$}}}}
\put(8.2,3.65){\makebox(0,0)[b]{{\mbox{\scriptsize $Sq^2$}}}}
\put(9.55,3.65){\makebox(0,0)[bl]{{\mbox{\scriptsize $Sq^2$}}}}
\put(11.1,3.65){\makebox(0,0)[bl]{{\mbox{\scriptsize $Sq^2$}}}}
\put(12,4.5){\makebox(0,0)[l]{{\mbox{\scriptsize 
              $H^{2n + 2}(X;\mbox{\bf Z$_2$})$}}}}
\put(12,3){\makebox(0,0)[l]{{\mbox{\scriptsize 
              $H^{2n}(X;\mbox{\bf Z$_2$})$}}}}
\put(12,1.5){\makebox(0,0)[l]{{\mbox{\scriptsize 
              $H^{2n - 2}(X;\mbox{\bf Z$_2$})$}}}}
\put(.3,1.3){\dashbox{.1}(2.9,.4){}}
\put(.3,2.8){\dashbox{.1}(2.9,.4){}}
\put(4.3,2.8){\dashbox{.15}(2.9,.4){}}
\put(8.3,2.8){\dashbox{.15}(2.9,.4){}}
\put(8.3,4.3){\dashbox{.15}(2.9,.4){}}
\put(1.7,.87){\makebox(0,0)[b]{{\mbox{\scriptsize $B_{2n - 2}$}}}}
\put(1.7,3.45){\makebox(0,0)[b]{{\mbox{\scriptsize $C_{2n}$}}}}
\put(5.7,2.4){\makebox(0,0)[b]{{\mbox{\scriptsize $D_{2n}$}}}}
\put(9.7,2.4){\makebox(0,0)[b]{{\mbox{\scriptsize $B_{2n}$}}}}
\put(9.7,4.95){\makebox(0,0)[b]{{\mbox{\scriptsize $C_{2n + 2}$}}}}
\end{picture}
\skp{.2}
\nd We conclude that the ring structure of $\ml$ determines the \a1
module structure. {\em Notice that the\/} \a1 {\em module structure of \/}
$H^{*}(X;\mbox{\bf Z$_2$})$ {\em can depend only on the map\/}
$\lambda \bmod{\,2}$.  
%
%
\skp{.2}
\newpage
\nd {\bf Example\/}. Let $P^3$ be the three dimensional cube and 
the map 
$$\lambda : {\cal F} \rightarrow {\mbox{{\bf Z}}}^3$$
($\bmod{\,2}$), be as in the diagram below.
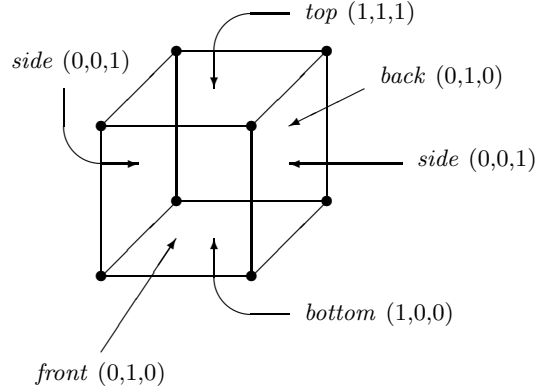
\skp{.2}
\begin{picture}(14,6)(0,0)
\put(4,2){\line(0,1){2}}
\put(4,2){\line(1,0){2}}
\put(6,2){\line(0,1){2}}
\put(4,4){\line(1,0){2}}
\put(4,2){\circle*{.15}}
\put(6,4){\circle*{.15}}
\put(6,2){\circle*{.15}}
\put(4,4){\circle*{.15}}
\put(5,3){\line(0,1){2}}
\put(5,3){\line(1,0){2}}
\put(7,3){\line(0,1){2}}
\put(5,5){\line(1,0){2}}
\put(5,3){\circle*{.15}}
\put(7,5){\circle*{.15}}
\put(7,3){\circle*{.15}}
\put(5,5){\circle*{.15}}
\put(4,2){\line(1,1){1}}
\put(4,4){\line(1,1){1}}
\put(6,2){\line(1,1){1}}
\put(6,4){\line(1,1){1}}

\put(5.5,2){\vector(0,1){.5}}
\put(6,2){\oval(1,1)[bl]}
\put(6.5,1.5){\line(-1,0){.5}}
\put(6.7,1.5){\makebox(0,0)[l]{{\mbox{\footnotesize {\em bottom\/} (1,0,0)}}}}

\put(5.5,5){\vector(0,-1){.5}}
\put(6,5){\oval(1,1)[tl]}
\put(6.5,5.5){\line(-1,0){.5}}
\put(6.7,5.5){\makebox(0,0)[l]{{\mbox{\footnotesize {\em top\/} (1,1,1)}}}}

\put(4,3.5){\vector(1,0){.5}}
\put(4,4){\oval(1,1)[bl]}
\put(3.5,4.5){\line(0,-1){.5}}
\put(3.6,4.85){\makebox(0,0){{\mbox{\footnotesize {\em side\/} (0,0,1)}}}}

\put(4,1){\vector(2,3){1}}
\put(4,.7){\makebox(0,0){{\mbox{\footnotesize {\em front\/} (0,1,0)}}}}

\put(7.5,4.5){\vector(-2,-1){1}}
\put(7.7,4.5){\makebox(0,0)[bl]{{\mbox{\footnotesize {\em back\/} (0,1,0)}}}}

\put(8,3.5){\vector(-3,0){1.5}}
\put(8.2,3.55){\makebox(0,0)[l]{{\mbox{\footnotesize {\em side\/} (0,0,1)}}}}

\end{picture}
\skp{0.2}
\nd Now
$$H^{*}(M^{6}(\lambda); \mbox{{\bf Z}}_2) = 
{\mbox{{\bf Z}}[v_1,v_2,\ldots,v_6]} \bigg/(I + J) \quad \bmod{\,2}$$
For $P^3$ we have $f_0 = 6$, $f_1 = 12$ and $f_2 = 8$ from which it 
follows easily that $h_0 = 1$, $h_1 = 3$, $h_2 = 3$ and $h_3 = 1$
where $h_i$ is the rank of 
$H^{2i}(M^{6}(\lambda); \mbox{{\bf Z}}_2)$.
The simplicial complex $K$ dual to $P^{3}$ is an octohedron with
vertices $\{v_1,v_2, \ldots, v_6\}$. The ideal of relations $I$ is 
generated by 
$v_{1}v_{6} = 0$, $v_{2}v_{4} = 0$ and $v_{3}v_{5} = 0$.
The ideal of relations $J$ is determined by the matrix representation 
$$\lambda = \pmatrix{1&0&0&0&0&1\cr
                     1&0&1&0&1&0\cr
                     1&1&0&1&0&0\cr}$$
This gives $v_{1} = v_{6} = v_{3} + v_{5} = v_{2} + v_{4}$.
Choose  as generators of 
$H^{2}(M^{6}(\lambda); \mbox{{\bf Z}}_2)$, $\{v_{1}, v_{2}, v_{3}\}$.
The relations in
$H^{4}(M^{6}(\lambda); \mbox{{\bf Z}}_2)$ become
${v_{1}}^{2} = 0$, ${v_{2}}^{2} = v_{1}v_{2}$ and 
${v_{3}}^{2} = v_{1}v_{3}$.
In $H^{6}(M^{6}(\lambda); \mbox{{\bf Z}}_2)$ we have 
$v_{1}{v_{2}}^2 = v_{1}{v_{3}}^2 = v_{3}{v_{1}}^2 = {v_{3}}^3 
= {v_{2}}^3 = 0$ and $v_{3}{v_{2}}^2 = v_{2}{v_{3}}^2 = v_{1}v_{2}v_{3}$.
We conclude that as \a1 modules
$$H^{*}(M^{6}(\ld); \mbox{\bf Z$_{2}$}) \; \cong \;
{\mathop{\oplus} \limits_{j = 0}^{3}}\; {\textstyle \sum^{2j}{S^0}} 
\oplus \;
2\,{\textstyle \sum^{2}{{\cal M}}}$$ 
In the next section we show that this 
is sufficient to enable us to read off $KO_{*}(M^{6}(\ld))$.
\skp{0.3}
\nd {\bf Problem.\/} Given \pn\/ and $\lambda$, find an algorithm which 
will determine the $Sq^2$ connections directly from the matrix 
representing $\lambda$, that is, {\em without\/} doing the algebra 
involved in solving the relations.
\skp{1}
\centerline{{\bf 4. The Adams Spectral Sequence for 
 $ko$-Homology}}
\skp{0.1}
\nd Let $X$ be {\em any\/} space with $H^*(X;\mbox{\bf Z$_2$})$ 
concentrated in even degrees. The ($\bmod{\,2}$) Adams spectral sequence 
relevant for our calculation takes the form
$$E_{2} \cong \mbox{Ext}_{\cal A}^{s,t}(H^{*}(ko \wedge X), 
\mbox{\bf Z$_{2}$})
\cong \mbox{Ext}_{{\cal A}(1)}^{s,t}(H^{*}(X), \mbox{\bf Z$_{2}$})
\Longrightarrow ko_{t - s}X$$
More details about this Adams spectral sequence can be found in,
for example, \cite{bb}.
\skp{0.2}
\nd At odd primes, in the case $X = M^{n}(\ld)$, the Atiyah-Hirzebruch 
spectral sequence converging to $ko_{*}X$ collapses for dimensional
reasons and we can conclude easily that $ko_{*}X$ has no odd torsion.
In fact, 
$$ko_{*}(M^{n}(\ld)) \otimes \mbox{\bf Z$_{(p)}$} \; \cong \;
H_{*}(M^{n}(\ld); \mbox{\bf Z$_{(p)}$}) \otimes ko_{*}$$
where {\bf Z$_{(p)}$} denotes the integers localized at $p$ odd.
So, a $\bmod{\,2}$ calculation suffices for the whole $ko$-theory.
\skp{0.2}
\nd Lemma 5 tells us that as \a1 modules
$$H^{*}(X; \mbox{\bf Z$_{2}$}) \; \cong \;
{\mathop{\oplus} \limits_{j = 0}^{k}}\; m_j\,{\textstyle \sum^{2j}{S^0}} 
\bigoplus \;
{\mathop{\oplus} \limits_{j = 0}^{l}}\; n_j\,{\textstyle 
\sum^{2j}{{\cal M}}}$$
where positive integers $m_j$ and $n_j$ denote the number of copies
of each summand located in dimension $2j$. Then
$$\mbox{Ext}_{{\cal A}(1)}^{s,t}(H^{*}(X), \mbox{\bf Z$_{2}$}) \; 
\cong \; {\mathop{\oplus} \limits_{j = 0}^{k}}\; 
m_j\! \cdot \mbox{Ext}_{{\cal A}(1)}^{s,t}({\textstyle \sum^{2j}{S^0}} , 
\mbox{\bf Z$_{2}$})
\bigoplus \;
{\mathop{\oplus} \limits_{j = 0}^{l}}\; 
n_j\! \cdot \mbox{Ext}_{{\cal A}(1)}^{s,t}({\textstyle 
\sum^{2j}{{\cal M}}} ,\mbox{\bf Z$_{2}$})$$
where the isomorphism is as 
$\mbox{Ext}_{{\cal A}(1)}^{s,t}(S^0, \mbox{\bf Z$_{2}$})$ modules.
\skp{0.2}
\nd The bigraded algebra
$\mbox{Ext}_{{\cal A}(1)}^{s,t}(S^0, \mbox{\bf Z$_{2}$})$ 
is well known, \cite{al}.
$$\mbox{Ext}_{{\cal A}(1)}^{s,t}(S^0, \mbox{\bf Z$_{2}$})\; \cong \;
\mbox{{\bf Z$_2$}}[a_0, a_1, w, b]\mbox{{\bf $\big/$}}
(a_{0}a_{1}, a_{1}^{3}, a_{1}w, w^2 + a_{0}^{2}b)$$
with $|a_0| = (0,1)$, $|a_1| = (1,1)$, $|w| = (4,3)$ and
$|b| = (8,4)$, where $|x| = (t - s,s)$ specifies the geometric degree
$t - s$ and the Adams filtration $s$. It's most easily represented
by the picture following. The vertical line segments indicate 
multiplication by $a_0$ and the sloping line segments, multiplication 
by $a_1$.
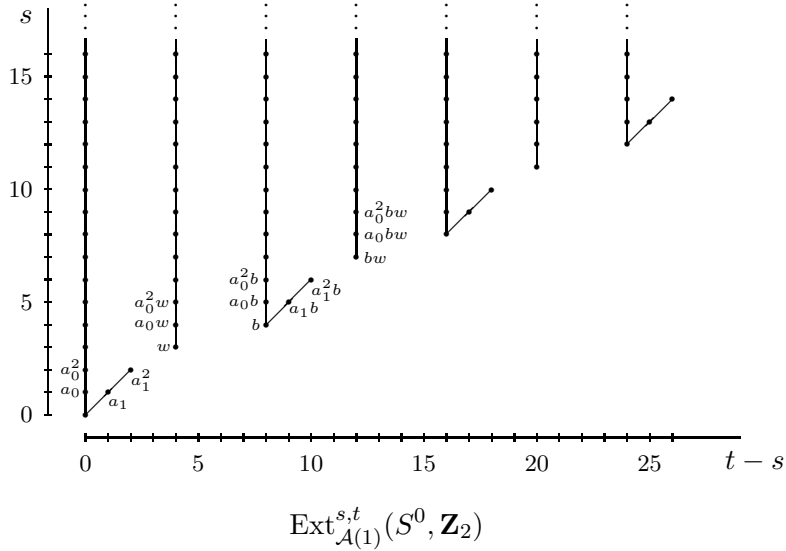
\skp{.1}
\newenvironment{ext_pt}{\begin{picture}(11,8)(-2,0)
\put(1.3,1){\line(1,0){8.7}}
\put(.8,1.3){\line(0,1){5.4}}
\multiput(1.3,.95)(.3,0){27}{\line(0,1){.1}}
\multiput(.75,1.3)(0,.3){17}{\line(1,0){.1}}
\put(1.3,1.3){\line(0,1){5}}
\put(.96,1.6){\makebox(0,0)[l]{{\mbox{\tiny $a_{0}$}}}}
\put(.96,1.9){\makebox(0,0)[l]{{\mbox{\tiny $a_{0}^{2}$}}}}
\multiput(1.3,6.45)(0,.15){3}{\circle*{.03}}
\multiput(1.3,1.3)(0,.3){17}{\circle*{.08}}
\put(1.6,1.51){\makebox(0,0)[tl]{{\mbox{\tiny $a_{1}$}}}}
\put(1.9,1.9){\makebox(0,0)[tl]{{\mbox{\tiny $a_{1}^{2}$}}}}

\put(1.3,1.3){\line(1,1){.6}}
\multiput(1.6,1.6)(.3,.3){2}{\circle*{.08}}

\put(2.5,2.2){\line(0,1){4.1}}
\put(2.243,2.2){\makebox(0,0)[l]{{\mbox{\tiny $w$}}}}
\put(1.95,2.5){\makebox(0,0)[l]{{\mbox{\tiny $a_{0}w$}}}}
\put(1.95,2.8){\makebox(0,0)[l]{{\mbox{\tiny $a_{0}^{2}w$}}}}

\multiput(2.5,6.45)(0,.15){3}{\circle*{.03}}
\multiput(2.5,2.2)(0,.3){14}{\circle*{.08}}

\put(3.7,2.5){\line(0,1){3.8}}
\multiput(3.7,6.45)(0,.15){3}{\circle*{.03}}     
\multiput(3.7,2.5)(0,.3){13}{\circle*{.08}}
\put(3.62,2.5){\makebox(0,0)[r]{{\mbox{\tiny $b$}}}}
\put(3.6,2.8){\makebox(0,0)[r]{{\mbox{\tiny $a_{0}b$}}}}
\put(3.6,3.1){\makebox(0,0)[r]{{\mbox{\tiny $a_{0}^{2}b$}}}}
\put(4,2.8){\makebox(0,0)[tl]{{\mbox{\tiny $a_{1}b$}}}}
\put(4.3,3.1){\makebox(0,0)[tl]{{\mbox{\tiny $a_{1}^{2}b$}}}}

\put(3.7,2.5){\line(1,1){.6}}
\multiput(4,2.8)(.3,.3){2}{\circle*{.08}}

\put(4.9,3.4){\line(0,1){2.9}}
\multiput(4.9,6.45)(0,.15){3}{\circle*{.03}}
\multiput(4.9,3.4)(0,.3){10}{\circle*{.08}}
\put(5,3.4){\makebox(0,0)[l]{{\mbox{\tiny $bw$}}}}
\put(5,3.7){\makebox(0,0)[l]{{\mbox{\tiny $a_{0}bw$}}}}
\put(5,4){\makebox(0,0)[l]{{\mbox{\tiny $a_{0}^{2}bw$}}}}

\put(6.1,3.7){\line(0,1){2.6}}
\multiput(6.1,6.45)(0,.15){3}{\circle*{.03}}
\multiput(6.1,3.7)(0,.3){9}{\circle*{.08}}

\put(6.1,3.7){\line(1,1){.6}}
\multiput(6.4,4)(.3,.3){2}{\circle*{.08}}

\put(7.3,4.6){\line(0,1){1.7}}
\multiput(7.3,6.45)(0,.15){3}{\circle*{.03}}
\multiput(7.3,4.6)(0,.3){6}{\circle*{.08}}

\put(8.5,4.9){\line(0,1){1.4}}
\multiput(8.5,6.45)(0,.15){3}{\circle*{.03}}
\multiput(8.5,4.9)(0,.3){5}{\circle*{.08}}

\put(8.5,4.9){\line(1,1){.6}}
\multiput(8.8,5.2)(.3,.3){2}{\circle*{.08}}

\setcounter{dim}{-5}
\multiput(1.3,.55)(1.5,0){6}{\addtocounter{dim}{5}\makebox(0,0)[b]
{{\footnotesize \arabic{dim}}}}
\put(10.2,.55){\makebox(0,0)[b]{$t - s$}}

\setcounter{dim}{-5}
\put(.6,6.6){\makebox(0,0)[r]{$s$}}
\multiput(.6,1.3)(0,1.5){4}{\addtocounter{dim}{5}\makebox(0,0)[r]
{\footnotesize {\arabic{dim}}}}}{\end{picture}}
\nd \begin{ext_pt}
\end{ext_pt}
\skp{-0.2}
\centerline{${\displaystyle
\mbox{Ext}_{{\cal A}(1)}^{s,t}(S^0, \mbox{\bf Z$_{2}$})}$}
\skp{0.5}
\nd The vertical multiplication by $a_0$ yields 
{\em multiplication-by-two\/}
extensions at $E_{\infty}$. The vertical towers in this diagram
produce copies of {\bf Z}$_{(2)}$, the integers localized at 2, in 
$ko_{*}S^0$. The other classes yield
copies of \z2. The class $b$ represents the Bott periodicity 
operator. Embedded in this picture then is $ko_{*}$ the coefficients
of $ko$-theory.
$$ko_{*}S^0 \cong \mbox{\zl2}
\; \oplus \; {\textstyle \sum^{1}}\,{\mbox{\z2}} 
\; \oplus \; {\textstyle \sum^{2}}\,{\mbox{\z2}} 
\; \oplus \; {\textstyle \sum^{4}}\,{\mbox{\zl2}}
\; \oplus \; {\textstyle \sum^{8}}\,{\mbox{\zl2}}
\; \oplus \; {\textstyle \sum^{9}}\,{\mbox{\z2}}
\; \oplus \; \ldots$$ 
$\mbox{Ext}_{{\cal A}(1)}^{s,t}({\cal M}, \mbox{\bf Z$_{2}$})$ 
is computed easily from 
$\mbox{Ext}_{{\cal A}(1)}^{s,t}({S^0}, \mbox{\bf Z$_{2}$})$ 
and the cofibration sequence associated to ${\cal M}$.
As a module over 
$\mbox{Ext}_{{\cal A}(1)}^{s,t}(S^0, \mbox{\bf Z$_{2}$})$, 
$\mbox{Ext}_{{\cal A}(1)}^{s,t}({\cal M}, \mbox{\bf Z$_{2}$})$ has 
generators $x,y,z,u$ with $|x| = (0,0)$, $|y| = (2,1)$, 
$|z| = (4,2)$ and  $|u| = (6,3)$ and relations 
$$a_{1}x = a_{1}y = a_{1}z = a_{1}u = 0,\;
a_{0}z = wx,\; a_{0}u = wy,\; wz = a_{0}bx, \;wu = a_{0}by$$

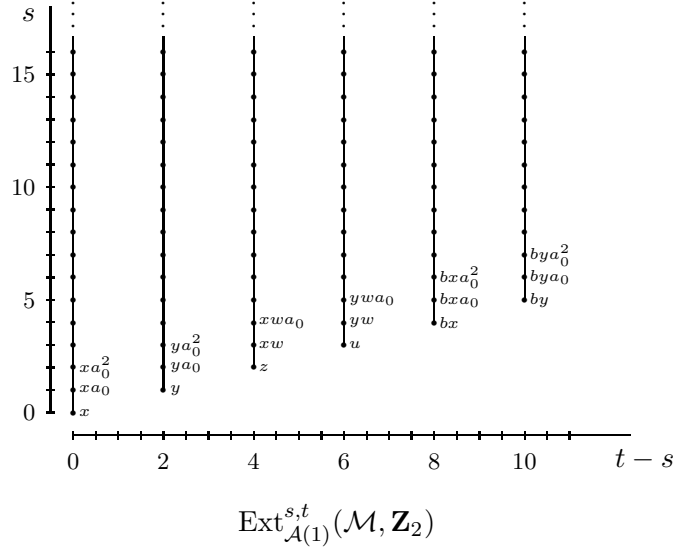
\skp{.6}

\newenvironment{ext_cp2}{\begin{picture}(10,7)(-2.5,0)


\put(1.3,1){\line(1,0){7.4}}
\put(1,1.3){\line(0,1){5.4}}
\multiput(1.3,.95)(.3,0){23}{\line(0,1){.1}}
\multiput(.95,1.3)(0,.3){17}{\line(1,0){.1}}

\put(1.3,1.3){\line(0,1){5}}
\multiput(1.3,6.45)(0,.15){3}{\circle*{.03}}
\multiput(1.3,1.3)(0,.3){17}{\circle*{.08}}
\put(1.38,1.3){\makebox(0,0)[l]{{\mbox{\tiny $x$}}}}
\put(1.38,1.6){\makebox(0,0)[l]{{\mbox{\tiny $xa_{0}$}}}}
\put(1.38,1.9){\makebox(0,0)[l]{{\mbox{\tiny $xa_{0}^{2}$}}}}

\put(2.5,1.6){\line(0,1){4.7}}
\multiput(2.5,6.45)(0,.15){3}{\circle*{.03}}
\multiput(2.5,1.6)(0,.3){16}{\circle*{.08}}
\put(2.59,1.6){\makebox(0,0)[l]{{\mbox{\tiny $y$}}}}
\put(2.59,1.9){\makebox(0,0)[l]{{\mbox{\tiny $ya_{0}$}}}}
\put(2.59,2.2){\makebox(0,0)[l]{{\mbox{\tiny $ya_{0}^{2}$}}}}

\put(3.7,1.9){\line(0,1){4.4}}
\multiput(3.7,6.45)(0,.15){3}{\circle*{.03}}
\multiput(3.7,1.9)(0,.3){15}{\circle*{.08}}
\put(3.77,1.9){\makebox(0,0)[l]{{\mbox{\tiny $z$}}}}
\put(3.77,2.2){\makebox(0,0)[l]{{\mbox{\tiny $xw$}}}}
\put(3.77,2.5){\makebox(0,0)[l]{{\mbox{\tiny $xwa_{0}$}}}}

\put(4.9,2.2){\line(0,1){4.1}}
\multiput(4.9,6.45)(0,.15){3}{\circle*{.03}}
\multiput(4.9,2.2)(0,.3){14}{\circle*{.08}}
\put(4.97,2.2){\makebox(0,0)[l]{{\mbox{\tiny $u$}}}}
\put(4.97,2.5){\makebox(0,0)[l]{{\mbox{\tiny $yw$}}}}
\put(4.97,2.8){\makebox(0,0)[l]{{\mbox{\tiny $ywa_{0}$}}}}

\put(6.1,2.5){\line(0,1){3.8}}
\multiput(6.1,6.45)(0,.15){3}{\circle*{.03}}
\multiput(6.1,2.5)(0,.3){13}{\circle*{.08}}
\put(6.17,2.5){\makebox(0,0)[l]{{\mbox{\tiny $bx$}}}}
\put(6.17,2.8){\makebox(0,0)[l]{{\mbox{\tiny $bxa_{0}$}}}}
\put(6.17,3.1){\makebox(0,0)[l]{{\mbox{\tiny $bxa_{0}^{2}$}}}}

\put(7.3,2.8){\line(0,1){3.5}}
\multiput(7.3,6.45)(0,.15){3}{\circle*{.03}}
\multiput(7.3,2.8)(0,.3){12}{\circle*{.08}}
\put(7.37,2.8){\makebox(0,0)[l]{{\mbox{\tiny $by$}}}}
\put(7.37,3.1){\makebox(0,0)[l]{{\mbox{\tiny $bya_{0}$}}}}
\put(7.37,3.4){\makebox(0,0)[l]{{\mbox{\tiny $bya_{0}^{2}$}}}}

\setcounter{dim}{-2}
\multiput(1.3,.55)(1.2,0){6}{\addtocounter{dim}{2}\makebox(0,0)[b]
{\footnotesize {\arabic{dim}}}}
\put(8.9,.55){\makebox(0,0)[b]{$t - s$}}

\setcounter{dim}{-5}
\put(.8,6.6){\makebox(0,0)[r]{$s$}}
\multiput(.8,1.3)(0,1.5){4}{\addtocounter{dim}{5}\makebox(0,0)[r]
{\footnotesize {\arabic{dim}}}}}{\end{picture}}

\nd \begin{ext_cp2}
\end{ext_cp2}
\skp{-0.2}
\centerline{${\displaystyle
\mbox{Ext}_{{\cal A}(1)}^{s,t}({\cal M}, \mbox{\bf Z$_{2}$})}$}
\skp{0.5}
\nd Since $\sum^{2}{{\cal M}} \simeq H^{*}(CP^2, \mbox{\bf Z$_{2}$})$ 
and noting that no differentials 
are possible in the spectral sequence, we can read off the connective
$ko$-homology of the complex projective plane
$$ko_{*}CP^2 \cong {\textstyle \sum^{2}}\,{\mbox{\zl2}} 
\; \oplus \; {\textstyle \sum^{4}}\,{\mbox{\zl2}} 
\; \oplus \; {\textstyle \sum^{6}}\,{\mbox{\zl2}}
\; \oplus \; {\textstyle \sum^{8}}\,{\mbox{\zl2}}
\; \oplus \; \ldots$$ 
%
The decomposition above of 
$\mbox{Ext}_{{\cal A}(1)}^{s,t}(H^{*}(X), \mbox{\bf Z$_{2}$})$
implies that its diagram
is obtained by superimposing shifted copies of the diagrams
for
$\mbox{Ext}_{{\cal A}(1)}^{s,t}(S^0, \mbox{\bf Z$_{2}$})$
and 
$\mbox{Ext}_{{\cal A}(1)}^{s,t}({\cal M}, \mbox{\bf Z$_{2}$})$.
Dimensional considerations and the fact that $d_r$ is a derivation
with respect to the action of 
$\mbox{Ext}_{{\cal A}(1)}^{s,t}(S^0, \mbox{\bf Z$_{2}$})$
allow us to conclude that one type of non-zero differential 
$$d_r : E_{r}^{s,t} \longrightarrow  E_{r}^{s + r, t + r -1}$$
is possible in the spectral sequence. It occurs on a copy of
$\mbox{Ext}_{{\cal A}(1)}^{s,t}(S^0, \mbox{\bf Z$_{2}$})$ as in
the diagram below. In the diagram we have identified the generator
$$c_{2j} \in 
\mbox{Ext}_{{\cal A}(1)}^{0,2j}(H^{*}(X), \mbox{\bf Z$_{2}$})$$
of an $\mbox{Ext}_{{\cal A}(1)}^{s,t}(\sum^{2j}{S^0}, \mbox{\bf Z$_{2}$})$ 
summand, with the dual of $c_{2j} \in C_{2j} \subseteq 
H^{2j}(X;\mbox{\bf Z$_2$})$. The class $\tilde{c}_{2p}$ represents some
linear combination of classes in 
$\mbox{Ext}_{{\cal A}(1)}^{0,2p}(H^{*}(X), \mbox{\bf Z$_{2}$})$
\skp{.2}
\newenvironment{ext_diff}{\begin{picture}(8,8)(-2,0)

\put(3.1,1){\line(1,0){3.9}}
\put(1,2.0){\line(0,1){1.6}}
\put(1,4.5){\line(0,1){2}}
\multiput(3.3,.95)(.3,0){12}{\line(0,1){.1}}
\multiput(.95,2.2)(0,.3){5}{\line(1,0){.1}}
\multiput(.95,4.7)(0,.3){5}{\line(1,0){.1}}
\multiput(1,3.92)(0,.13){3}{\circle*{.03}}

\put(4.2,4.9){\line(0,1){1.4}}
\multiput(4.2,6.45)(0,.13){3}{\circle*{.03}}
\multiput(4.2,4.9)(0,.3){5}{\circle*{.09}}

\put(4.2,4.9){\line(1,1){.6}}
\multiput(4.2,4.9)(.3,.3){3}{\circle*{.09}}
\put(4.05,4.9){\makebox(0,0)[r]{{\mbox{\scriptsize 
$b^{j}\tilde{c}_{2p}$}}}}

\put(4.8,2.2){\line(0,1){1.4}}
\multiput(4.8,3.72)(0,.13){3}{\circle*{.03}}
\put(5,3.95){\makebox(0,0)[l]{{\mbox{\footnotesize $d_{r}$}}}}
\multiput(4.8,2.2)(0,.3){5}{\circle*{.09}}

\put(4.8,2.2){\line(1,1){.6}}
\multiput(4.8,2.2)(.3,.3){3}{\circle*{.09}}
\put(4.65,2.2){\makebox(0,0)[r]{{\mbox{\scriptsize $b^{k}c_{2q}$}}}}

\put(7.5,.5){\makebox(0,0)[b]{$t - s$}}
\put(.8,6.6){\makebox(0,0)[r]{$s$}}

\put(4.5,5.1){\vector(-1,4){.01}}
\put(4.8,5.4){\vector(-1,4){.01}}

\multiput(4.8,2.2)(-.01,.1){30}{\circle*{.01}}
\multiput(5.1,2.5)(-.01,.1){30}{\circle*{.01}}}{\end{picture}}

\nd \begin{ext_diff}
\end{ext_diff}
\skp{-0.2}
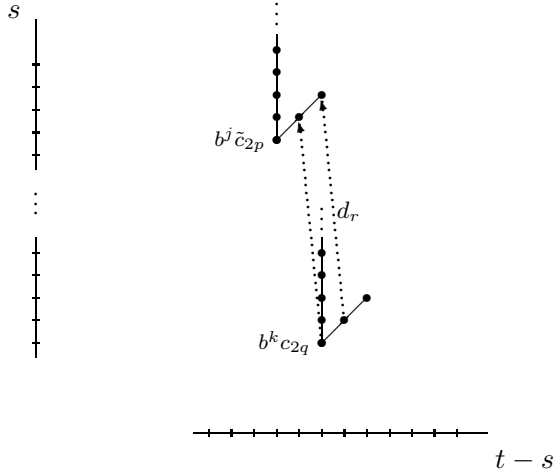
\centerline{A Differential in the Adams Spectral Sequence for
${\displaystyle ko_{*}X}$}
\skp{0.5}
\nd {\bf Important Remark.\/} Since $b$ has $(t - s, s)$ bidegree
$(8,4)$, this differential {\em cannot\/}
occur in the Adams Spectral Sequence for a toric manifold or
toric variety, of dimension less than 12, with $\bmod{\,2}$ 
cohomology concentrated in even degrees. Consequently, 
the spectral sequence collapses without any further analysis
and theorem 1 holds for such spaces.
\skp{0.2}
\nd We shall use the fact that a toric manifold is a manifold to 
prove that there can be no non-zero differentials in the spectral
sequence. Choose $q$ minimal so that for some $r$,
we have
$d_{r}(b^{k}c_{2q}) \neq 0$. Next, choose the smallest such $r$
so that for some $k$, 
we have $d_{r}(b^{k}c_{2q}) \neq 0$. The derivation 
property  of $d_{r}$ with respect to multiplication by the
periodicity operator $b$, implies then that $d_{r}(c_{2q}) \neq 0$ 
and so we can assume that $k = 0$.
%
%
%
%
%
%
%
%
%
%
%

We restrict now to the case $X = \ml$ a toric manifold of dimension 
$2n$. Consider {\em all\/} $2q$ dimensional submanifolds $M_{F_{i}}$
of $\ml$ corresponding to $q$-faces $F_{i}$. The inclusions 
$$M_{F_{i}} \hookrightarrow \ml$$
induce  maps of Adams Spectral Sequences and in particular, maps
$$\mbox{Ext}_{{\cal A}(1)}^{s,t}(H^{*}(M_{F_{i}}), \mbox{\bf Z$_{2}$})
\longrightarrow
\mbox{Ext}_{{\cal A}(1)}^{s,t}(H^{*}(\ml), \mbox{\bf Z$_{2}$})$$
In each 
$\mbox{Ext}_{{\cal A}(1)}^{0,2q}(H^{*}(M_{F_{i}}), \mbox{\bf Z$_{2}$})$
there is a unique class corresponding to the fundamental class
$[M_{F_{i}}]$. Theorem 3 tells us that $c_{2q}$ is a linear combination
of the images of the classes $[M_{F_{i}}]$. Because 
$d_{r}(c_{2q}) \neq 0$, the naturality of the Adams Spectral Sequence
implies that
$d_{r}([M_{F_{i}}]) \neq 0$ for some $i$. In other words, a $q$-face
$F = F_{i}$ of $P^{n}$ must exist with a non-zero differential in the
Adams Spectral Sequence for $ko_{*}(M_F)$ supported on the top class of
filtration zero. We shall use the result following to show that this 
cannot be the case for the manifold $M_F$ and so complete the proof of
theorem 1 
\skp{0.4}
\nd {\bf Theorem 6.\/} {\em Let \/} $M$ {\em be an orientable manifold 
of dimension\/} $n$ 
{\em Then\/} $M$ {\em is a spin manifold if the top 
dimensional cohomology class is not in the image of\/} $Sq^2$.
\skp{0.25}
\nd {\bf Proof:\/} Let $v \in H^{*}(M)$ be the total Wu class of $M$.
It satisfies the property that $Sq(v) = w$ where $Sq$ is the total
Steenrod operation and $w$ is the total Stiefel-Whitney class. Since
$M$ is orientable
we have $v_2 = w_2$ where
$w_2$ is the second Stiefel-Whitney class. The Wu formula for $M$,
(\cite{hu}, page 261), is
$$<a \cup v, [M]> \; = \; <Sq(a), [M]>$$
for any $a \in H^{*}(M)$. In particular, for any class 
$x \in H^{n - 2}(M)$, we have
$$<x \cup w_{2}, [M]> \; = \; <x \cup v_{2}, [M]> \; = \; <Sq^{2}(x), [M]>$$
So, if $Sq^{2}(x) = 0$ for all $x$ we must have $w_2 = 0$ by 
Poincar\'{e} duality and so $M$ is a spin manifold. \block
\skp{0.4}
\nd {\bf Corollary 7.\/} {\em There are no non-zero differentials in the 
Adams Spectral Sequence for\/} $ko_{*}(M_F)$ {\em supported on the top class
in filtration zero\/}.
\skp{0.25}
\nd {\bf Proof:\/} Suppose such a differential did exist. Then
the ${\cal A}(1)$ module $H^{*}(M_{F}), \mbox{\bf Z$_{2}$})$
must contain a summand $S^{0}$ in the top dimension $2q$. In 
particular, the top class in $H^{2q}((M_{F}), \mbox{\bf Z$_{2}$})$
is not in the image of $Sq^{2}$ and so $M_F$ must be spin manifold.
This implies, (\cite{abp}), that $M_F$ is {\em orientable\/} with 
respect to $ko_{*}$. We can now apply Poincar\'{e}-Lefschetz duality,
(\cite{st}, page 39(a)), to conclude that as a $ko_{*}$ module,
$ko_{*}(M_F)$ must contain a summand, free on a single generator
in $ko_{2q}(M_F)$ dual to the single summand on the generator in
$ko^{0}(M_F)$. This contradicts the existence of the differential.
\block
\skp{0.4}
\nd The fact that the Adams spectral sequence collapses leaves us 
with possible group extension problems before we can read off
the group $ko_{*}(M^{n}(\ld))$. Fortunately, in our case these are
not difficult. As mentioned earlier, the vertical multiplication
by $a_0$ yields {\em multiplication-by-two\/} extensions at 
$E_{\infty}$. All other classes in the spectral sequence are products
of $a_1$. Vertical extensions across copies of $ko_{*}(S^0)$,
of {\bf Z}$_2$ groups to groups of higher torsion, cannot occur 
because products of $a_1$ yield elements of order two in $ko$-theory. 
\skp{0.4}
\nd We conclude that, if as \a1 modules
$$H^{*}(M^{n}(\ld); \mbox{\bf Z$_{2}$}) \; \cong \;
{\mathop{\oplus} \limits_{j = 0}^{k}}\; m_j\,{\textstyle \sum^{2j}{S^0}} 
\bigoplus \;
{\mathop{\oplus} \limits_{j = 0}^{l}}\; n_j\,{\textstyle \sum^{2j}{M}}$$
then 
$$ko_{*}(M^{n}(\ld)) \; \cong \;
{\mathop{\oplus} \limits_{j = 0}^{k}}\; 
m_j\,{\textstyle \sum^{2j}{ko_{*}S^0}} 
\bigoplus \;
{\mathop{\oplus} \limits_{j = 0}^{l}}\; 
n_j\,{\textstyle \sum^{2j}{ko_{*}M}}$$
where the graded groups $ko_{*}{S^0}$ and $ko_{*}M$ are described above.

Our calculation shows that multiplication by the Bott element $b$ is a
monomorphism in $E_{\infty}$ and hence in $ko_{*}(M^{n}(\ld))$. So, we can
invert $b$ to get the periodic $KO$-homology of $M^{n}(\ld)$.
$$KO_{*}(M^{n}(\ld)) \; \cong \;
{\mathop{\oplus} \limits_{j = 0}^{k}}\; 
m_j\,{\textstyle \sum^{2j}{KO_{*}S^0}} 
\bigoplus \;
{\mathop{\oplus} \limits_{j = 0}^{l}}\; 
n_j\,{\textstyle \sum^{2j}{KO_{*}M}}$$
where
$$KO_{*}S^0 \;\; \cong \;\; \ldots  
\: \oplus \: {\textstyle \sum^{-6}}\,{\mbox{\z2}}
\: \oplus \: {\textstyle \sum^{-4}}\,{\mbox{{\bf Z}}}
\: \oplus \: \mbox{{\bf Z}}
\: \oplus \: {\textstyle \sum^{1}}\,{\mbox{\z2}} 
\: \oplus \: {\textstyle \sum^{2}}\,{\mbox{\z2}} 
\: \oplus \: {\textstyle \sum^{4}}\,{\mbox{{\bf Z}}}
\: \oplus \: \ldots$$ 
and 
$$KO_{*}M \;\; \cong \;\; \ldots  
\: \oplus \: {\textstyle \sum^{-4}}\,{\mbox{{\bf Z}}}
\: \oplus \: {\textstyle \sum^{-2}}\,{\mbox{{\bf Z}}}
\: \oplus \: \mbox{{\bf Z}}
\: \oplus \: {\textstyle \sum^{2}}\,{\mbox{{\bf Z}}}
\: \oplus \: {\textstyle \sum^{4}}\,{\mbox{{\bf Z}}}
\: \oplus \: {\textstyle \sum^{6}}\,{\mbox{{\bf Z}}}
\: \oplus \: \ldots$$ 
\newpage
\centerline{{\bf 5. The $KO$-cohomology of Toric Manifolds}}
\skp{0.1}
\nd We employ the universal coefficient exact sequence following to
compute the $KO$-cohomology from the $KO$-homology.
\skp{.4}
\nd {\bf Theorem 8. [D. W. Anderson, \cite{and}, theorem 2.4]} 
{\em Let $X$ be a CW-complex. For all $n$,
there is a natural exact sequence\/}
$$0 \; \rightarrow \; \mbox{lim}^{1}KO^{m -1}(X) \; \rightarrow \;
\mbox{Ext}_{{\mbox{\bf Z}}}(KSp_{m - 1}(X), \mbox{\bf Z})
\; \rightarrow \quad$$
$$\quad \quad \quad \quad \mbox{lim}^{0}KO^{m}(X) \; \rightarrow \;
\mbox{Hom}_{{\mbox{\bf Z}}}(KSp_{m}(X), \mbox{\bf Z})
\; \rightarrow \; 0$$
{\em where these limits are over the filtration of $X$ by finite 
subcomplexes\/}. \block
\skp{0.2}
\nd In our case, $X = M^{n}(\ld)$ is a finite complex and we 
are left with the sequence
$$0 \; \rightarrow \;
\mbox{Ext}_{{\mbox{\bf Z}}}(KSp_{m - 1}M^{n}(\ld), \mbox{\bf Z})
\; \rightarrow \; KO^{m}M^{n}(\ld) \; \rightarrow \;
\mbox{Hom}_{{\mbox{\bf Z}}}(KSp_{m}M^{n}(\ld), \mbox{\bf Z})
\; \rightarrow \; 0$$
\nd Bott periodicity implies 
$KSp_{m}M^{n}(\ld) \; \cong \; KO_{m - 4}M^{n}(\ld)$. Combining this
with the results of the previous section, namely, that the groups 
$KO_{*}M^{n}(\ld)$ are direct sums of copies of {\bf Z} and \z2,
we see that the short exact sequence splits. Explicitly, if
$KO_{m}M^{n}(\ld) \; \cong \; \alpha_{m} \cdot \mbox{{\bf Z}}
\; \oplus \; \beta_{m} \cdot$ \z2, for integers $\alpha_m$
and $\beta_m$, then, as groups
$$KO^{m}M^{n}(\ld) \; \cong \; \alpha_{m - 4} \cdot \mbox{\bf Z} \;
\oplus \; \beta_{m - 5} \cdot \mbox{\z2}$$
\skp{.2}
We conclude with a remark about the module structure. Let $DM^{n}(\ld)$ 
denotes the $S$-dual of $M^{n}(\ld)$. If
$$H^{*}(M^{n}(\ld); \mbox{\bf Z$_{2}$}) \; \cong \;
{\mathop{\oplus} \limits_{j = 0}^{k}}\; m_j\,{\textstyle \sum^{2j}{S^0}} 
\bigoplus \;
{\mathop{\oplus} \limits_{j = 0}^{l}}\; n_j\,{\textstyle \sum^{2j}{M}}$$
then by duality
$$H^{*}(DM^{n}(\ld); \mbox{\bf Z$_{2}$}) \; \cong \;
{\mathop{\oplus} \limits_{j = 0}^{k}}\; m_j\,{\textstyle \sum^{- 2j}{S^0}} 
\bigoplus \;
{\mathop{\oplus} \limits_{j = 0}^{l}}\; n_j\,{\textstyle \sum^{2j - 2}{M}}$$
So, except for dimension shifts. the Adams spectral sequence for 
$ko_{*}DM^{n}(\ld)$ looks much as it did for $ko_{*}M^{n}(\ld)$
We cannot use the same arguments however to conclude that the spectral
sequence collapses. Instead, we now know the groups $KO^{m}M^{n}(\ld)$ 
and so we can use a rank argument to conclude that all differentials must 
be zero. This allows us to read off $ko_{*}DM^{n}(\ld)$ as a $ko_{*}S^{0}$
module because we know the $ko_{*}S^{0}$  module structure of $ko_{*}M$. 
Again, the Bott element $b$ acts as a monomorphism and we can conclude the 
$KO_{*}S^0$ module structure of $KO_{*}DM^{n}(\ld)$ and so of
$KO^{*}M^{n}(\ld)$


\skp{0.6}
\nd {\small {\sc Department of Mathematics, Rider University, Lawrenceville,
New Jersey, 08648}}
\skp{0.1}
\nd {\small {\em E-mail address\/}: bahri@rider.edu}
\skp{0.4}
\nd {\small {\sc Department of Mathematics, Hunter College, New York, 
New York 10021}}
\skp{0.1}
\nd {\small {\em E-mail address\/}: mbenders@shiva.hunter.cuny.edu}
\end{document}